\documentclass[12pt]{article}
\newcommand{\qed}{\hfill\rule{4pt}{8pt}\par\vspace{\baselineskip}}
\newtheorem{de}{Definition}[section]

\newtheorem{pr}[de]{Proposition}
\newtheorem{co}[de]{Corollary}
\newtheorem{re}[de]{Remark}

\newtheorem{te}[de]{Theorem}

\def\bea{\begin{eqnarray*}}\def\eea{\end{eqnarray*}}

\begin{document}
\title{Doubles of (quasi) Hopf algebras \\
and some examples of quantum groupoids and vertex groups related to them}

\author { Florin Panaite\\Institute of Mathematics of the Romanian Academy\\
P. O. Box 1-764, RO-70700 Bucharest, Romania\\e-mail: fpanaite@stoilow.imar.ro
}
\date{}
\maketitle
\section{Introduction}
${\;\;\;}$Let $A$ be a finite dimensional Hopf algebra, 
$D(A)=A^{*cop}\otimes A$ its Drinfel'd double and $\cal H$$(A)=A\#A^*$ its 
Heisenberg double. The relation between $D(A)$ and $\cal H$$(A)$ has been 
found by J.-H. Lu in \cite{lu1} (see also \cite{mon}, p. 196): the  
multiplication of $\cal H$$(A)$ may be obtained by twisting the multiplication 
of $D(A)^*$ by a certain left 2-cocycle which in turn is obtained from 
the $R$-matrix of $D(A)$. It was also obtained in \cite{lu1} that 
$\cal H$$(A)$ becomes a left $D(A)$-module algebra under a certain action 
of $D(A)$ on $\cal H$$(A)$ (formula (35) in \cite{lu1}). \\
${\;\;\;}$All these may be obtained alternatively using a more direct 
approach, which also shows that the above mentioned action of $D(A)$ on 
$\cal H$$(A)$ is manifestly the left re\-gu\-lar action of $D(A)$ on $D(A)^*$  
(by identifying $\cal H$$(A)$ and $D(A)^*$ as linear spaces). 
The general setting is the following:  
if $(H, R)$ is a quasitriangular bialgebra and we define a new  
multiplication on $H^*$ by $f\cdot g=\sum (R^2\rightharpoonup g)
(R^1\rightharpoonup f)$, then $H^*$ with this new multiplication (denoted in 
what follows by $H^*_R$) becomes a left $H$-module algebra under the 
left regular action $\rightharpoonup $ (this is well-known, see also  
\cite{ae}, \cite{maj3}, \cite{gz}, \cite{bpvo} for some more general 
versions in terms of Drinfel'd twists). Although very simple, this 
construction may have some nice applications, for instance $H^*_R$ may be 
noncommutative even if $H$ was cocommutative$-$it was discovered recently in   
\cite{w1}, \cite{w2} that an important algebra arising in  
noncommutative string theory is an example of this type; this 
dis\-co\-ve\-ry has  
also been applied to noncommutative quantum field theory in \cite{o}. And,   
if $A$ is a finite dimensional Hopf algebra and $H=D(A)$, then $H^*_R$ is 
just $\cal H$$(A)$. \\
${\;\;\;}$For reasons to be discussed below, we were not satisfied with the 
description of $H^*_R$ as a left $H$-module algebra and we were led to 
consider also the right regular action of $H$ on $H^*_R$. It turns out that 
$H^*_R$ is a right $H^{cop}$-module algebra (that is, the right action 
satisfies a ``reversed Leibniz rule''), so $H^*_R$ is an algebra in the 
tensor category of $H-H^{cop}$-bimodules (we say that it is an 
$H-H^{cop}$-bimodule algebra). If we endow this category with the  
braiding given by multiplying to the left by $R_{21}$ (as usual) and from the 
right by $R^{-1}$, we shall prove that $H^*_R$ is quantum commutative as an 
algebra in this braided tensor category (which is also equivalent to saying 
that $H^*_R$ is a quantum commutative left $H\otimes H^{op\;cop}$-module  
algebra). In particular $\cal H$$(A)$ is  
a quantum commutative $D(A)-D(A)^{cop}$-bimodule algebra and we like to think  
of this as the most natural ``tensor categorical'' interpretation of the 
Heisenberg double.\\
${\;\;\;}$In section 4 we discuss some more facts about $H^*_R$, Drinfel'd 
doubles and Heisenberg doubles. For instance, we discuss the relation 
between the multiplication of $H^*_R$ and Majid's ``covariantised product'' 
and, for a finite dimensional Hopf algebra $A$, we give a formula for the 
canonical element $W\in $$\cal H$$(A)\otimes $$\cal H$$(A)$, solution to 
the pentagon equation, in terms of the $R$-matrix of $D(A)$ and the map 
$Q$ expressing the factorizability of $D(A)$ (the formula is: 
$W=(Q^{-1}\otimes Q^{-1})(R_{21})$; unfortunately this formula does not seem 
to offer an answer to the following natural question: is there an 
``explanation'', in terms of the structure of the Drinfel'd double only, for 
the fact that $W$ is a solution to the pentagon equation on the Heisenberg 
double?). \\
${\;\;\;}$In section 5 we speak about quantum groupoids. The general concepts 
of bialgebroid and Hopf algebroid (=quantum groupoid) have been 
introduced by J.-H. Lu in \cite{lu2}, with inspiration and motivation coming 
from Poisson geometry and by generalizing previous ones (\cite{rav}, 
\cite{mal}) where the base algebra was assumed to be commutative (we refer 
to \cite{en}, \cite{nv}, \cite{sch}, \cite{bm} for discussions concerning 
the relation  
between these concepts and other objects known as ``quantum groupoids'', 
such as weak Hopf algebras \cite{bns} and Takeuchi's $\times _R$-bialgebras  
\cite{tak}). In \cite{lu2} Lu proved that if $A$ is a finite dimensional 
Hopf algebra and $V$ is a quantum commutative left $D(A)$-module algebra then 
$V\# A$ is a Hopf algebroid over $V$, and that $A^*$ is a quantum commutative 
left $D(A)$-module algebra, so that $\cal H$$(A^*)$ is a quantum groupoid 
over $A^*$. We would like to obtain a quantum groupoid having $\cal H$$(A)$ 
as $base$, and we proceed as follows: first we generalize Lu's theorem, by 
proving that if $(H,R)$ is any quasitriangular Hopf algebra and $V$ is a 
quantum commutative left $H$-module algebra then $V\#R_{(l)}$ is a quantum  
groupoid over $V$, where $R_{(l)}$ is a certain finite dimensional Hopf 
subalgebra of $H$ (Radford's notation). Then, if $A$ is a finite dimensional 
Hopf algebra, $\cal H$$(A)$ is a quantum commutative 
$D(A)-D(A)^{cop}$-bimodule algebra, hence it is a quantum commutative 
left $D(A)\otimes D(A)^{op\;cop}$-module algebra so that the above result 
may be applied and we obtain a Hopf algebroid with $\cal H$$(A)$ as base 
and $\cal H$$(A)\#(A\otimes A^{*op})$ as total algebra. \\
${\;\;\;}$Let us mention that a generalization of Lu's theorem has been 
independently obtained also very recently in \cite{bm}.\\
${\;\;\;}$In section 6 we give a slight generalization of the concept of 
``vertex group'' introduced by Richard Borcherds in his recent (partly 
Hopf-algebraic) approach to vertex algebras (see \cite{b1}, \cite{b2}), 
by allowing the ``ring of  
singular functions'' to be noncommutative. With this terminology, we prove 
that if $A$ is a finite dimensional cocommutative Hopf algebra, then 
$\cal H$$(A)$ has some properties making it a vertex group over $A$. These 
properties are natural from the Hopf-algebraic point of view, but we do not  
know whether this example may be relevant for the theory of vertex algebras. \\
${\;\;\;}$Now, we come back to our starting point, namely the construction 
of $H^*_R$. We have tried to perform it for quasi-bialgebras $H$ instead of 
bialgebras, the multiplication on $H^*$ being defined by the same formula. 
Naturally, this multiplication is $not$ associative in general, 
but surprisingly  
$H^*_R$ is also $not$, in general, an algebra in the tensor category of left  
$H$-modules, as for  
bialgebras. We have tried to find a tensor category in which $H^*_R$ lives 
as an algebra, and we found the category of $H-H^{cop}$-bimodules. So, in 
section 2, the construction and properties of $H^*_R$ are given directly 
for quasi-bialgebras. This greater generality not only indicates what is the 
most natural tensor-categorical interpretation for the case of bialgebras, 
but, in view of the fact that $\cal H$$(A)=D(A)^*_R$ for Hopf algebras, 
suggests a possible definition for the Heisenberg double of a finite 
dimensional quasi-Hopf algebra $A$, as $D(A)^*_R$, where $D(A)$ is the 
quantum double of $A$ introduced in \cite{maj2}, \cite{hn1}, \cite{hn2}. 
In section 7 we compute explicitly this $D(A)^*_R$ for the case when 
$D(A)$ is a slight generalization of the Dijkgraaf-Pasquier-Roche 
quasi-Hopf algebra $D^{\omega }(G)$ introduced in \cite{dpr}.\\ 
${\;\;\;}$This  
definition of the Heisenberg double as a $non-associative$ algebra (but 
which is an algebra in a certain tensor category) may seem rather strange,  
so let us mention that such non-associative algebras occur naturally 
in the literature on quasi-Hopf algebras, in various contexts such as 
smash products (\cite{bpvo}), cohomology and deformation theory 
(\cite{ms}, \cite{ss}))  
and algebraic quantum field theory (\cite{mack}). A somehow dual situation 
appears in \cite{am1}, \cite{am2}, where it was proposed, as a general 
philosophy, to try to study non-associative algebras by expressing them,  
when possible, as algebras in certain tensor categories (especially ones 
associated to quasi-Hopf algebras). For instance, the octonions and higher 
Cayley algebras may be studied in this framework, see \cite{am1}.   
\section{Preliminaries}
${\;\;\;}$In this section, we recall some 
definitions and fix the notation that 
will be used in the rest of the paper. Throughout, $k$ will
be a fixed field and all algebras, linear spaces etc.
will be over $k$; unadorned $\otimes $ means
$\otimes _k$. For coalgebras and Hopf algebras, we shall use
the framework of \cite{sw2}; in particular, for coalgebras,
we shall use $\Sigma -$notation: $\Delta (h)=\sum h_1\otimes
h_2$, $(I\otimes \Delta )(\Delta (h))=(\Delta \otimes I)
(\Delta (h))=\sum h_1\otimes h_2\otimes h_3$, etc.
\begin{de} (\cite{dr2}) Let $H$ be a $k-$algebra, $\Delta: H\rightarrow  
H\otimes H$, $\varepsilon: H\rightarrow k$ two algebra 
homomorphisms. $H$ is called a quasi-bialgebra if there exists 
an invertible element $\Phi \in H\otimes H\otimes H$ such that, 
for all elements $h\in H$, we have:
$$(I\otimes \Delta)(\Delta (h))=\Phi ((\Delta \otimes I)
(\Delta (h))\Phi ^{-1}$$
$$(\varepsilon \otimes I)(\Delta (h))=h \;and\; 
(I\otimes \varepsilon )(\Delta (h))=h$$
$$(I\otimes I\otimes \Delta)(\Phi )(\Delta \otimes I\otimes I)
(\Phi )=(1\otimes \Phi )(I\otimes \Delta\otimes I)(\Phi )
(\Phi \otimes 1)$$
$$(I\otimes \varepsilon \otimes I)(\Phi )=1\otimes 1$$
where $I=id_H$. The map $\Delta $ is called the coproduct or the 
comultiplication, $\varepsilon $ the counit and $\Phi $ the associator. \\
${\;\;\;}$$H$ is called a quasi-Hopf algebra if, moreover,
there exist an anti-automorphism $S$ of the algebra $H$  
and elements $\alpha $ and $\beta $ of $H$ such that, for all 
$h\in H$, we have: 
$$\sum S(h_1)\alpha h_2=\varepsilon (h)\alpha \; and \; 
\sum h_1\beta S(h_2)=\varepsilon (h)\beta $$
$$\sum X^1\beta S(X^2)\alpha X^3=1 \; and \; 
\sum S(x^1)\alpha x^2\beta S(x^3)=1$$
where $\Phi =\sum X^1\otimes X^2\otimes X^3$,
$\Phi ^{-1}=\sum x^1\otimes x^2\otimes x^3$ (formal notation), 
and we used also the $\Sigma$-notation :
$\Delta (h)=\sum h_1\otimes h_2$. In this case, $S$ is called 
the antipode of $H$.
\end{de}
${\;\;\;}$Let us note that every Hopf algebra with bijective antipode
is a quasi-Hopf algebra with $\Phi =1\otimes 1\otimes 1$ and 
$\alpha =\beta =1$.\\
${\;\;\;}$We note the following two consequences of the definitions of
$S, \alpha ,\beta $:
$\varepsilon (\alpha )\varepsilon (\beta )=1$, $\varepsilon \circ
S=\varepsilon $. Moreover, the axioms imply that 
$(\varepsilon \otimes I\otimes I)(\Phi )=
(I\otimes I\otimes \varepsilon )(\Phi )=1$.\\
\begin{de} (\cite{dr2}) A quasi-bialgebra or a quasi-Hopf algebra 
$H$ is termed quasitriangular if there exists an invertible element 
$R\in H\otimes H$ such that: 
$$(\Delta \otimes I)(R)=\Phi _{312}R_{13}\Phi ^{-1}_{132}R_{23}\Phi $$
$$(I\otimes \Delta)(R)=\Phi^{-1}_{231}R_{13}\Phi_{213}R_{12}\Phi
^{-1}$$
$$\Delta ^{cop}(h)=R\Delta (h)R^{-1}\; for \;all \; h\in H$$ 
where, if $t$ denotes a permutation of $\{1, 2, 3\}$, then we set
$\Phi _{t(1)t(2)t(3)}=\sum X^{t^{-1}(1)}\otimes X^{t^{-1}(2)}\otimes
X^{t^{-1}(3)}$
and $R_{ij}$ means $R$ acting non-trivially in the $i^{th}$ and $j^{th}$
positions of $H\otimes H\otimes H$. 
\end{de}
${\;\;\;}$If $R$ satisfies these conditions it is 
called an $R-$matrix. From these relations one can deduce 
the quasi-Yang-Baxter equation:
$$R_{12}\Phi_{312}R_{13}\Phi ^{-1}_{132}R_{23}\Phi =\Phi _{321}R_{23}\Phi
^{-1}_{231}R_{13}\Phi _{213}R_{12}$$
${\;\;\;}$Also, it is easy to see that:
$$(\varepsilon \otimes I)(R)=(I\otimes \varepsilon )(R)=1$$
${\;\;\;}$As a general rule, the tensor components of the associator 
$\Phi $ of a quasi-bialgebra will be denoted using big letters, for instance 
$$\Phi =\sum X^1\otimes X^2\otimes X^3=\sum Y^1\otimes Y^2\otimes Y^3 \;etc$$
and the ones of $\Phi ^{-1}$ with small letters, for instance 
$$\Phi ^{-1}=\sum x^1\otimes x^2\otimes x^3=\sum y^1\otimes y^2\otimes y^3 
\;etc$$
${\;\;\;}$An $R$-matrix will be usually denoted by 
$$R=\sum R^1\otimes R^2=\sum r^1\otimes r^2 $$
${\;\;\;}$If $(H, R)$ is a quasitriangular Hopf algebra and $B$ is a left 
$H$-module algebra (i.e. an algebra in the tensor category  
$H-mod$) then $B$ is called $quantum \;commutative$ in \cite{cw} if it 
is commutative as algebra in the braided tensor category of left $H$-modules, 
i.e. if $bb'=\sum (R^2\cdot b')(R^1\cdot b)$ for all $b, b'\in B$. We 
extend this terminology and we call $quantum \;commutative$ any algebra $B$ 
in a braided tensor category which is commutative with respect to the 
braiding $c$ of the category, namely $m_B\circ c_{B,B}=m_B$ where $m_B$ is 
the multiplication of $B$.   
\section{The main result for quasi-bialgebras}
${\;\;\;}$Let $H$ be a quasi-bialgebra and denote by $H_l$, $H_r$ and $H_{lr}$ 
the categories of left $H$-modules, right $H$-modules and $H$-bimodules 
respectively. In these categories we introduce tensor products, as follows. 
If $V,W\in H_l$ then $V\otimes W\in H_l$ with $h\cdot (v\otimes w)=
\Delta (h)\cdot (v\otimes w)=\sum h_1\cdot v\otimes h_2\cdot w$. If 
$V,W\in H_r$ then $V\otimes W\in H_r$ with $(v\otimes w)\cdot h=(v\otimes w)
\cdot \Delta ^{cop}(h)=\sum v\cdot h_2\otimes w\cdot h_1$. If $V,W\in H_{lr}$ 
then $V\otimes W\in H_{lr}$ with $h\cdot (v\otimes w)\cdot h'=\Delta (h)
\cdot (v\otimes w)\cdot \Delta ^{cop}(h')=\sum h_1\cdot v\cdot h'_2\otimes 
h_2\cdot w\cdot h'_1$.\\
${\;\;\;}$It is well-known (see \cite{kass}, Chapter XV) that $H_l$ with 
the above tensor product becomes a tensor category, with associativity 
constraints given by 
$$\Phi _{U,V,W}:(U\otimes V)\otimes W\rightarrow U\otimes (V\otimes W)$$
$$\Phi _{U,V,W}((u\otimes v)\otimes w)=\Phi \cdot (u\otimes (v\otimes w))=
\sum X^1\cdot u\otimes (X^2\cdot v\otimes X^3\cdot w)$$
for all $U,V,W\in H_l$. Similarly one can prove that $H_r$ and $H_{lr}$ 
become also tensor categories, with associativity constraints given by
$$\Phi _{U,V,W}((u\otimes v)\otimes w)=(u\otimes (v\otimes w))\cdot 
\Phi _{321}=\sum u\cdot Y^3\otimes (v\cdot Y^2\otimes w\cdot Y^1)$$
$$\Phi _{U,V,W}((u\otimes v)\otimes w)=\Phi \cdot (u\otimes (v\otimes w))
\cdot \Phi _{321}=\sum X^1\cdot u\cdot Y^3\otimes (X^2\cdot v\cdot Y^2  
\otimes X^3\cdot w\cdot Y^1)$$
for $U,V,W\in H_r$ and $U,V,W\in H_{lr}$ respectively (the unit constraints  
are the usual ones). \\
${\;\;\;}$Suppose now that $(H,R)$ is a quasitriangular quasi-bialgebra. 
It is well-known (see \cite{kass}) that the tensor category $H_l$ is 
braided, the braiding being given by 
$$c_{V,W}:V\otimes W\rightarrow W\otimes V$$
$$c_{V,W}(v\otimes w)=\sum R^2\cdot w\otimes R^1\cdot v$$
for $V,W\in H_l$. Similarly one can prove that $H_r$ and $H_{lr}$ become 
also braided tensor categories, the braidings being given by
$$c_{V,W}(v\otimes w)=\sum w\cdot U^1\otimes v\cdot U^2$$
$$c_{V,W}(v\otimes w)=\sum R^2\cdot w\cdot U^1\otimes R^1\cdot v\cdot U^2$$
for $V,W\in H_r$ and $V,W\in H_{lr}$ respectively, where $U=\sum U^1\otimes 
U^2$ is the inverse of $R$. \\
${\;\;\;}$Suppose again that $(H,R)$ is a quasitriangular quasi-bialgebra 
and consider the left and right regular actions of $H$ on $H^*$, that is 
$(h\rightharpoonup p)(h')=p(h'h)$, $(p\leftharpoonup h)(h')=p(hh')$, 
for $p\in H^*$ and $h,h'\in H$. Obviously $H^*$ is an $H$-bimodule with 
these actions.\\
${\;\;\;}$On $H^*$ we can consider the convolution product, given by 
$(fg)(h)=\sum f(h_1)g(h_2)$ (which is $not$ associative in general). We  
introduce another product in $H^*$, given by 
$$f\cdot g=\sum (R^2\rightharpoonup g)(R^1\rightharpoonup f)$$
which is also $not$ associative in general. Denote by $H^*_R$ the pair 
$(H^*, \cdot)$. Then we have the following   
\begin{te}
a) $H^*_R$ is an algebra in the tensor category $H_{lr}$ (we shall say that 
it is an $H-H^{cop}$-bimodule algebra), that is, for all  
$f,g,l\in H^*$ and $h,h'\in H$ we have: 
$$h\rightharpoonup (f\cdot g)\leftharpoonup h'=\sum (h_1\rightharpoonup f
\leftharpoonup h'_2)\cdot (h_2\rightharpoonup g\leftharpoonup h'_1)$$
$$(f\cdot g)\cdot l=\sum (X^1\rightharpoonup f\leftharpoonup Y^3)\cdot  
((X^2\rightharpoonup g\leftharpoonup Y^2)\cdot (X^3\rightharpoonup l
\leftharpoonup Y^1))$$
$$\varepsilon \cdot f=f\cdot \varepsilon=f$$
$$h\rightharpoonup \varepsilon\leftharpoonup h'=\varepsilon (h)
\varepsilon (h')\varepsilon $$
b) $H^*_R$ is quantum commutative as an algebra in the braided tensor 
category $H_{lr}$, that is, for all $f,g\in H^*$, we have 
$$f\cdot g =\sum (R^2\rightharpoonup g\leftharpoonup U^1)\cdot 
(R^1\rightharpoonup f\leftharpoonup U^2)$$
\end{te}
{\bf Proof:} for $f,g\in H^*$, the product $f\cdot g$ is given by:
$$(f\cdot g)(h)=\sum g(h_1R^2)f(h_2R^1)$$
for all $h\in H$. We shall prove first $b)$. We calculate:\\[2mm]
$((R^2\rightharpoonup g\leftharpoonup U^1)\cdot (R^1\rightharpoonup f
\leftharpoonup U^2))(h)=\\[2mm]
=\sum (R^1\rightharpoonup f\leftharpoonup U^2)(h_1r^2)(R^2\rightharpoonup g
\leftharpoonup U^1)(h_2r^1)\\[2mm]
=\sum f(U^2h_1r^2R^1)g(U^1h_2r^1R^2)\\[2mm]
=\sum f(U^2r^2h_2R^1)g(U^1r^1h_1R^2)$\\[2mm]
(using the relation $\Delta ^{cop}(h)R=R\Delta (h)$)\\[2mm]
$=\sum f(h_2R^1)g(h_1R^2)\\[2mm]
=(f\cdot g)(h)$\\[2mm] 
${\;\;\;}$Now we prove $a)$. We have: \\[2mm]
$(h\rightharpoonup f\cdot g\leftharpoonup h')(h'')=\\[2mm]
=(f\cdot g)(h'h''h)\\[2mm]
=\sum g(h'_1h''_1h_1R^2)f(h'_2h''_2h_2R^1)$\\[3mm]
$(\sum (h_1\rightharpoonup f\leftharpoonup h'_2)\cdot (h_2\rightharpoonup g
\leftharpoonup h'_1))(h'')=\\[2mm]
=\sum (h_2\rightharpoonup g\leftharpoonup h'_1)(h''_1R^2)(h_1\rightharpoonup 
f\leftharpoonup h'_2)(h''_2R^1)\\[2mm]
=\sum g(h'_1h''_1R^2h_2)f(h'_2h''_2R^1h_1)\\[2mm]
=\sum g(h'_1h''_1h_1R^2)f(h'_2h''_2h_2R^1)$\\[2mm]
using again the relation $\Delta ^{cop}(h)R=R\Delta (h)$. \\
${\;\;\;}$For the second relation, we calculate:\\[2mm]
$((f\cdot g)\cdot l)(h)=\sum l(h_1R^2)(f\cdot g)(h_2R^1)=\\[2mm]
=\sum l(h_1R^2)g((h_2)_1(R^1)_1r^2)f((h_2)_2(R^1)_2r^1)\\[2mm]
=\sum l(h_1X^1R^2x^2\rho ^2Y^3)g((h_2)_1X^2R^1x^1Y^1r^2)f((h_2)_2X^3x^3
\rho ^1Y^2r^1)$\\[2mm]
(using the relation $\;(\Delta \otimes id)(R)=\sum X^2R^1x^1Y^1\otimes 
X^3x^3\rho ^1Y^2\otimes X^1R^2x^2\rho ^2Y^3$, where $R=\rho =\sum \rho ^1
\otimes \rho ^2$)\\[2mm]
$=\sum l(h_1X^1r^2x^2\rho ^2Y^3)g((h_2)_1R^2X^3x^3\rho ^1Y^2)f((h_2)_2R^1
X^2r^1x^1Y^1)$\\[2mm]
using the quasi-Yang-Baxter equation. We compute now the right hand side 
evaluated in $h$:\\[2mm]
$(\sum (X^1\rightharpoonup f\leftharpoonup Y^3)\cdot ((X^2\rightharpoonup g
\leftharpoonup Y^2)\cdot (X^3\rightharpoonup l\leftharpoonup Y^1)))(h)=\\[2mm]
=\sum ((X^2\rightharpoonup g\leftharpoonup Y^2)\cdot (X^3\rightharpoonup l
\leftharpoonup Y^1))(h_1R^2)(X^1\rightharpoonup f\leftharpoonup Y^3)
(h_2R^1)\\[2mm]
=\sum l(Y^1(h_1)_1(R^2)_1r^2X^3)g(Y^2(h_1)_2(R^2)_2r^1X^2)
f(Y^3h_2R^1X^1)\\[2mm]
=\sum l(h_1Y^1(R^2)_1r^2X^3)g((h_2)_1Y^2(R^2)_2r^1X^2)f((h_2)_2Y^3R^1X^1)$
\\[2mm]
(using the relation $(id \otimes \Delta )(\Delta (a))=\Phi (\Delta \otimes 
id)(\Delta (a))\Phi ^{-1}$)\\[2mm]
$=\sum l(h_1Y^1y^1T^1R^2x^2r^2X^3)g((h_2)_1Y^2y^2\rho ^2T^3x^3r^1X^2)
f((h_2)_2Y^3y^3\rho ^1T^2R^1x^1X^1)$\\[2mm]
(using the relation $\;(id \otimes \Delta )(R)=\sum y^3\rho ^1T^2R^1x^1
\otimes y^1T^1R^2x^2\otimes y^2\rho ^2T^3x^3$, where $\rho =R$, 
$\Phi =\sum T^1\otimes T^2\otimes T^3$)\\[2mm]
$=\sum l(h_1T^1R^2x^2r^2X^3)g((h_2)_1\rho ^2T^3x^3r^1X^2)f((h_2)_2\rho ^1T^2
R^1x^1X^1)$\\[2mm]
and this is obviously equal to the expression obtained for 
$((f\cdot g)\cdot l)(h)$. \\
${\;\;\;}$The relations $\varepsilon \cdot f=f\cdot \varepsilon =f$ and 
$h\rightharpoonup \varepsilon \leftharpoonup h'=\varepsilon (h)
\varepsilon (h')\varepsilon $ are obvious, using the fact that 
$\sum \varepsilon (h_1)h_2=h=\sum h_1\varepsilon (h_2)$ and 
$\sum \varepsilon (R^1)R^2=\sum R^1\varepsilon (R^2)=1$.\qed
\begin{re}{\em Although $H^*_R$ is an algebra in the tensor category 
$H_{lr}$,  
in general it is $not$ an algebra in the tensor categories $H_l$ or $H_r$.} 
\end{re}
\begin{re}{\em Suppose that $F:(H,R)\rightarrow (H',R')$ is a morphism of  
quasitriangular quasi-bialgebras and consider the transpose $F^*:H'^*_{R'}
\rightarrow H^*_R$. It is easy to see that the map $F^*$ has the 
following properties:
$$F^*(f\cdot g)=F^*(f)\cdot F^*(g)$$
$$F^*(F(h)\rightharpoonup f)=h\rightharpoonup F^*(f)$$
$$F^*(f\leftharpoonup F(h))=F^*(f)\leftharpoonup h$$
for all $f,g\in H'^*$ and $h\in H$. The map $F$ induces a tensor functor 
$F_*:H'_{lr}\rightarrow H_{lr}$ (see \cite{kass}) in the usual manner, 
and if $F$ is an isomorphism then this functor is an isomorphism of 
tensor categories and in this case the above relations imply that $H^*_R$ 
and $F_*(H'^*_{R'})$ are isomorphic as algebras in the tensor category 
$H_{lr}$. So, we can say that $H^*_R$ ``depends only on the isomorphism 
class of $(H,R)$''.}
\end{re} 
\section{The Heisenberg double $vs$ the Drinfel'd double for Hopf algebras}
${\;\;\;}$Let $(H, R)$ be a quasitriangular bialgebra. 
We can consider it as a  
quasitriangular quasi-bialgebra with trivial associator, so the result in 
the previous section may be applied to $H$. Since the associator of $H$ is 
trivial, $H^*_R$ is an associative algebra and is an algebra also in the 
tensor categories $H_l$ and $H_r$ (that is, a left $H$-module algebra and 
a right $H^{cop}$-module algebra), but in general is $not$ quantum  
commutative as algebra in $H_l$ or $H_r$. \\
${\;\;\;}$Let us note also the trivial fact that if $(H, R)$ is a 
quasitriangular bialgebra then $(H^{cop}, R_{21})$ is also a 
quasitriangular bialgebra and $(H^{cop})^*_{R_{21}}=(H^*_R)^{op}$.\\[2mm] 
${\;\;\;}$Let now $A$ be a finite dimensional Hopf algebra with antipode $S$. 
The Drinfel'd double of $A$, denoted by $D(A)$, is a quasitriangular 
Hopf algebra realized on the $k$-linear space $A^*\otimes A$; its coalgebra 
structure is the one of $A^{*cop}\otimes A$, the algebra structure is given by
$$(p\otimes a)(p'\otimes a')=\sum p(a_1\rightharpoonup p'\leftharpoonup 
S^{-1}(a_3))\otimes a_2a'$$
for all $p,p'\in A^*$ and $a,a'\in A$ (see \cite{rad1}), and the $R$-matrix 
is
$$R=\sum (\varepsilon \otimes e_i)\otimes (e^i\otimes 1)$$
where $\{e_i\}$ is a basis of $A$ and $\{e^i\}$ its dual basis in $A^*$. \\
${\;\;\;}$The Heisenberg double of $A$, denoted by $\cal H$$(A)$,  
is the smash product $A\# A^*$, where $A^*$ acts on $A$ via the 
left regular action $p\rightharpoonup a=\sum p(a_2)a_1$ for all $p\in A^*, 
a\in A$, so its multiplication is
$$(a\otimes p)(a'\otimes p')=\sum a(p_1\rightharpoonup a')\otimes p_2p'$$
for all $p,p'\in A^*$ and $a,a'\in A$, where $\Delta (p)=\sum p_1\otimes p_2$ 
is the coalgebra structure of $A^*$. The present name of $A\#A^*$ is from 
\cite{sts}, while in \cite{nill} and \cite{maj1} it was used under the 
heading ``Weyl algebra of $A$''. \\
${\;\;\;}$We have now the following application of the above construction 
$H^*_R$:
\begin{pr} If $A$ is a finite dimensional Hopf algebra and $H=D(A)$, then 
$H^*_R=\cal H$$(A)$, so $\cal H$$(A)$ is an algebra in the tensor categories 
$D(A)_l$, $D(A)_r$, $D(A)_{lr}$ and is quantum commutative as an algebra in 
$D(A)_{lr}$. 
\end{pr}
{\bf Proof:} the fact that the multiplication of $\cal H$$(A)$ may be obtained 
from the one of $D(A)^*$ is due to J.-H. Lu in \cite{lu1}, but with a 
slightly different approach, so we include here a proof for completeness.\\ 
The multiplication in $D(A)^*_R$ is given by:
$$(b\otimes g)\cdot (a\otimes f)=\sum (R^2\rightharpoonup (a\otimes f))
(R^1\rightharpoonup (b\otimes g))$$
for all $a,b\in A$ and $f,g\in A^*$. Let $x\in A$, $p\in A^*$ and denote by 
$\langle ,\rangle $ the evaluation map; we calculate:\\[2mm]
$\langle (b\otimes g)\cdot (a\otimes f), p\otimes x\rangle =\\[2mm]
=\sum \langle a\otimes f, (p\otimes x)_1R^2\rangle \langle b\otimes g, 
(p\otimes x)_2R^1\rangle \\[2mm]
=\sum \langle a\otimes f, (p_2\otimes x_1)(e^i\otimes 1)\rangle \langle 
b\otimes g, (p_1\otimes x_2)(\varepsilon \otimes e_i)\rangle \\[2mm]
=\sum \langle a\otimes f, p_2(x_1\rightharpoonup e^i\leftharpoonup S^{-1}(x_3))
\otimes x_2\rangle \langle b\otimes g, p_1\otimes x_4e_i\rangle\\[2mm]
=\sum \langle p_2, a_1\rangle \langle e^i, S^{-1}(x_3)a_2x_1\rangle \langle f, 
x_2\rangle \langle p_1, b\rangle \langle g, x_4e_i\rangle \\[2mm]
=\sum \langle p, ba_1\rangle \langle f, x_2\rangle \langle g, x_4S^{-1}(x_3)
a_2x_1\rangle \\[2mm]
=\sum \langle p, ba_1\rangle \langle f, x_2\rangle \langle g, a_2x_1 
\rangle \\[2mm]
=\sum \langle p, ba_1\rangle \langle g_1, a_2\rangle \langle g_2, x_1\rangle 
\langle f, x_2\rangle \\[2mm]
=\sum \langle p, b(g_1\rightharpoonup a)\rangle \langle g_2f, x\rangle \\[2mm]
=\sum \langle b(g_1\rightharpoonup a)\otimes g_2f, p\otimes x \rangle $\\[2mm]
hence $(b\otimes g)\cdot (a\otimes f)=\sum b(g_1\rightharpoonup a)
\otimes g_2f$, q.e.d. \qed
\begin{re}{\em 
We would like to emphasize that the actions of $D(A)$ on $\cal H$$(A)$ are 
just the left and right regular actions of $D(A)$ on $\cal H$$(A)$ identified 
as vector spaces with $D(A)^*$. These actions may be described explicitly, 
using the following formula for the comultiplication of $D(A)^*$ (which may 
be found in \cite{lu1}):
$$\Delta _{D(A)^*}(a\otimes p)=\sum (a_1\otimes e^ip_1e^j)\otimes 
(S^{-1}(e_j)a_2e_i\otimes p_2)$$
where $\{e_i\}$ is a basis in $A$ and $\{e^i\}$ its dual basis in $A^*$. 
They look as follows:
$$(p\otimes b)\rightharpoonup (a\otimes q)=\sum p_2(a_2)q_2(b)(a_1\otimes 
p_3q_1S^{*-1}(p_1))$$
$$(a\otimes q)\leftharpoonup (p\otimes b)=\sum p(a_1)q_1(b_2)(S^{-1}(b_3)
a_2b_1\otimes q_2)$$
for all $p, q\in A^*$ and $a, b\in A$ (the first one is formula (35) in 
\cite{lu1}).} 
\end{re}
\begin{re}{\em 
The fact that $H^*_R$ is a left $H$-module algebra via the left regular 
action (for any quasitriangular bialgebra $(H,R)$) may be obtained 
alternatively using the more general framework of twisting the multiplication 
of a module algebra by a Drinfel'd twist (see \cite{ae}, \cite{maj3}, 
\cite{gz}, \cite{bpvo}). We use the Drinfel'd twist $(R_{21})^{-1}$ and from 
the above mentioned results it follows that $(H^*_R)^{op}$ is a left module 
algebra over $H^{cop}$ (hence $H^*_R$ is a left module algebra over $H$), 
and also using \cite{bpvo}, Prop. 2.17 we obtain that the algebras 
$H^*\#H$ and $(H^*_R)^{op}\#H^{cop}$ are isomorphic. In particular, for a 
finite dimensional Hopf algebra $A$, we obtain that the Heisenberg double of 
$D(A)^*$ (which is $D(A)^*\#D(A)$) is isomorphic to $\cal H$$(A)^{op}\#D(A)
^{cop}$ as algebras.\\
${\;\;\;}$ Let us mention that the Heisenberg double and the Drinfel'd 
double appear also together in \cite{cr} where they join to build the 
algebra $X=$$\cal H$$(A)\otimes D(A)^{op}$ having the nice property that 
there exists a $vector\;space-preserving$ equivalence of categories between 
the categories of Hopf bimodules over $A$ and of left $X$-modules.}
\end{re}
\begin{re}{\em 
Let $(H, R)$ be a quasitriangular bialgebra. Then the braided   
tensor category $H_r$ of right $H^{cop}$-modules may be identified with 
the braided tensor ca\-te\-go\-ry of left $H^{op\;cop}$-modules, where the  
$R$-matrix of $H^{op\;cop}$ is $(R_{21})^{-1}$. Moreover, if $V$ is an 
algebra, then $V$ is a right $H^{cop}$-module algebra if and only if $V$ is a 
left $H^{op\;cop}$-module algebra, and it is left quantum commutative if 
and only if it is right quantum commutative. \\
${\;\;\;}$Consequently, since $H^*_R$ is a quantum commutative 
$H-H^{cop}$-bimodule algebra, we obtain that $H^*_R$ is a quantum commutative 
left $H\otimes H^{op\;cop}$-module algebra, where the $R$-matrix of 
$H\otimes H^{op\;cop}$ is $\sum (R^1\otimes U^2)\otimes (R^2\otimes U^1)$ 
with $U=R^{-1}$ and the left action of $H\otimes H^{op\;cop}$ on $H^*_R$ is  
of course given by $(h\otimes h')\cdot p=h\rightharpoonup p\leftharpoonup h'$ 
for all $h, h'\in H$ and $p\in H^*$. 
In particular, if $A$ is a finite dimensional Hopf  
algebra, then $\cal H$$(A)$ is a quantum commutative left 
$D(A)\otimes D(A)^{op\;cop}$-module algebra. We shall use this in the 
next section.}
\end{re}   
\begin{pr} Let $(H,R)$ be a quasitriangular bialgebra. Then $H^*_R$ is 
quantum commutative as a left $H$-module algebra if and only if 
$R=1\otimes 1$. In particular, if $A$ is a finite dimensional Hopf algebra 
then $\cal H$$(A)$ is a quantum commutative left $D(A)$-module algebra 
if and only if $dim (A)=1$. 
\end{pr}
{\bf Proof:} if $f,g\in H^*$ then $f\cdot g=\sum (R^2\rightharpoonup g)
\cdot (R^1\rightharpoonup f)$ if and only if 
$$\sum g(h_1r^2)f(h_2r^1)=\sum f(h_1r^2R^1)g(h_2r^1R^2)$$
for all $h\in H$, and this holds for all $f,g\in H^*$ if and only if 
$$\sum h_2r^1\otimes h_1r^2=\sum h_1r^2R^1\otimes h_2r^1R^2$$
for all $h\in H$, which, by applying the relation 
$\Delta ^{cop}(h)R=R\Delta (h)$ several  
times is equivalent to 
$$\sum R^1h_1\otimes R^2h_2=\sum R^2r^1h_1\otimes R^1r^2h_2$$
for all $h\in H$, and this is of course equivalent to $R=1\otimes 1$. 
The second statement follows immediately from the first.\qed
\begin{re}{\em  
Let $A$ be a finite dimensional Hopf algebra. Then the left  
action of $D(A)$ on $\cal H$$(A)$ restricts to an action of $D(A)$ on 
$A^*\subseteq $$\cal H$$(A)$, given by 
$$(p\otimes b)\rightharpoonup q=\sum q_2(b)p_2q_1S^{*-1}(p_1)$$
So, $A^*$ is a left $D(A)$-module algebra under this action, and it was 
proved in \cite{lu2} that it is quantum commutative, a fact which 
allowed Lu to introduce a quantum groupoid structure on $\cal H$$(A^*)$ 
(we shall  
come back to the subject of quantum groupoids in the next section).\\
${\;\;\;}$Let us mention that $A$ is also a quantum commutative left 
$D(A)$-module algebra under the action: 
$$(p\otimes a)\cdot b=\sum (a_1bS(a_2))\leftharpoonup S^{-1}(p)$$
for all $a, b\in A$ and $p\in A^*$. This action is nicely obtained in 
\cite{z} by considering the right regular action of $D(A)$ on $D(A)^*$, 
restricting it to an action on $A^{op}$ and then transforming it via $S$ 
into a left  
action on $A$.}  
\end{re}
${\;\;\;}$If $(H,R)$ is a quasitriangular Hopf algebra, we shall denote by 
$Q$ the map $Q:H^*\rightarrow H$, $Q(p)=\sum p(R^2r^1)R^1r^2$. With this 
notation, let us recall the following
\begin{de} (\cite{rsts}) If $(H,R)$ is a finite dimensional 
quasitriangular Hopf algebra, it is called factorizable if the map $Q$ is 
a linear isomorphism. 
\end{de}
${\;\;\;}$This is a natural condition, which was also proved in \cite{h} 
to be useful in defining certain invariants for 3-manifolds. \\
${\;\;\;}$It is natural to see how far is the map $Q$ 
from being an algebra map, if we 
consider on $H^*$ the algebra structure $H^*_R$.    
\begin{pr} If $f,g\in H^*$, then 
$$Q(f\cdot g)=\sum Q(R^1\rightharpoonup f)Q(g\leftharpoonup R^2)$$
\end{pr}
{\bf Proof:} let us compute (we denote $R=r=\rho =T=\alpha =\beta =\gamma  
=\delta $):\\[2mm]
$Q(f\cdot g)=\sum (f\cdot g)(R^2r^1)R^1r^2\\[2mm]
=\sum g((R^2)_1(r^1)_1T^2)f((R^2)_2(r^1)_2T^1)R^1r^2\\[2mm]
=\sum f(\gamma ^2\beta ^1T^1)g(\delta ^2\alpha ^1T^2)\gamma ^1\delta ^1
\alpha ^2\beta ^2$\\[2mm]
(using: $(\Delta \otimes id)(R)=\sum \alpha ^1\otimes \beta ^1\otimes 
\alpha ^2\beta ^2$ and $(id\otimes \Delta )(R)=\sum \gamma ^1\delta ^1
\otimes \delta ^2\otimes \gamma ^2$)\\[2mm]
$=\sum f(\gamma ^2T^1\beta ^1)g(\delta ^2T^2\alpha ^1)\gamma ^1\delta ^1
\beta ^2\alpha ^2$\\[2mm]
(using the Yang-Baxter equation: $\sum \beta ^1T^1\otimes \alpha ^1T^2
\otimes \alpha ^2\beta ^2=\sum T^1\beta ^1\otimes T^2\alpha ^1\otimes 
\beta ^2\alpha ^2$)\\[2mm]
$=\sum f(\gamma ^2\beta ^1T^1)g(T^2\delta ^2\alpha ^1)\gamma ^1\beta ^2
\delta ^1\alpha ^2$\\[2mm]
(using the Yang-Baxter equation: $\sum T^1\beta ^1\otimes \delta ^1\beta ^2
\otimes \delta ^2T^2=\sum \beta ^1T^1\otimes \beta ^2\delta ^1\otimes T^2
\delta ^2$)\\[2mm]
$=\sum (T^1\rightharpoonup f)(\gamma ^2\beta ^1)\gamma ^1\beta ^2 
(g\leftharpoonup T^2)(\delta ^2\alpha ^1)\delta ^1\alpha ^2\\[2mm]
=\sum Q(T^1\rightharpoonup f)Q(g\leftharpoonup T^2)$, q.e.d. \qed
${\;\;\;}$It is natural to try to define an algebra structure on $H^*$ such 
that, with respect to this structure, $Q$ becomes an algebra map. In view 
of the previous proposition and the fact that $R^{-1}=\sum R^1\otimes 
S^{-1}(R^2)$, it is natural to define
$$f\underline{\cdot }g=\sum (R^1\rightharpoonup f)\cdot 
(g\leftharpoonup S(R^2))$$
and indeed one can now check that this multiplication is associative (due to 
the associativity of $\cdot $) and that $Q(f\underline{\cdot }g)=
Q(f)Q(g)$ (due to the previous proposition). This multiplication 
$\underline{\cdot }$ is just Majid's ``covariantised product'' (see 
\cite{maj1}, Th. 7.4.1), as one can easily check. The fact that 
$Q(f\underline{\cdot }g)=Q(f)Q(g)$ is also proved in \cite{maj1}, 
Prop. 7.4.3.\\
${\;\;\;}$Let us say few more words about the relation between the two 
multiplications $\cdot $ and $\underline{\cdot }$ on $H^*$. Define 
the following action of $H$ on $H^*$:
$$h\triangleright f=\sum h_2\rightharpoonup f\leftharpoonup S(h_1)$$
Obviously $H^*$ is a left $H$-module via this action; using the description 
of $\underline{\cdot }$ in terms of $\cdot $ and the fact that  
$(H^*_R, \cdot )$ is an algebra in the category of $H-H^{cop}$-bimodules, 
one can see immediately that $(H^*, \underline{\cdot })$ is a left 
$H$-module algebra with respect to the action $\triangleright $. 
It is $not$ quantum commutative, but satisfies a  
condition called by Majid ``braided commutativity'' (see \cite{maj1}, 
Exp. 9.4.10). It is nice that this condition may be easily proved using the 
fact that $(H^*_R, \cdot )$ is quantum commutative as an $H-H^{cop}$-bimodule  
algebra. \\
${\;\;\;}$Of course, if $H$ is finite dimensional and factorizable, we can 
express the covariantised product as
$$f\underline{\cdot }g=Q^{-1}(Q(f)Q(g))$$
${\;\;\;}$Now let $A$ be a finite dimensional Hopf algebra. It is well-known 
that the Drinfel'd double of $A$ is factorizable. An explicit proof (for the 
explicit realization of the double we work with) is written down in 
\cite{rad2}. It uses the following
\begin{pr} (\cite{rad2}) Let $(H,R)$ be a finite dimensional quasitriangular 
Hopf algebra, and denote by $u=\sum S(R^2)R^1$ the canonical Drinfel'd 
element of $H$. Then $(H,R)$ is factorizable if and only if 
$u^{-1}\leftharpoonup H^*=H$, or equivalently $H$ is a free right 
$H^*$-module with generator $u^{-1}$. Moreover, in this case, if $y\in H$, 
then $Q^{-1}(y)=p\leftharpoonup u^{-1}$, where $p\in H^*$ is the unique 
element satisfying the relation $u^{-1}y=u^{-1}\leftharpoonup p$. 
\end{pr}
${\;\;\;}$Another criterion for factorizability may be found in \cite{gw}.\\
${\;\;\;}$Actually, one can give a direct proof of the factorizability of 
$D(A)$: it is easy to see that in this case the map $Q$ is given by 
$$Q(a\otimes p)=\sum a_1\rightharpoonup p\leftharpoonup S^{-1}(a_3)
\otimes a_2$$
for all $a\in A, p\in A^*$, and it has an inverse given by
$$Q^{-1}:D(A)\rightarrow D(A)^*$$
$$Q^{-1}(p\otimes a)=\sum a_2\otimes S^{-1}(a_1)\rightharpoonup p
\leftharpoonup a_3$$     
We can also write down explicit formulae for Radford's criterion applied to 
a Drinfel'd double:
\begin{pr} Let $p\otimes a\in D(A)$; then we have
$$u^{-1}(p\otimes a)=u^{-1}\leftharpoonup (\sum S^{-2}(a_2)\otimes S^{-1}(a_1)
\rightharpoonup p)$$
Consequently, we have:
$$Q^{-1}(p\otimes a)=(\sum S^{-2}(a_2)\otimes S^{-1}(a_1)\rightharpoonup p)
\leftharpoonup u^{-1}$$
\end{pr}
{\bf Proof:} a direct computation, using the formula $u^{-1}=\sum e^i\otimes 
S^2(e_i)$.\qed 
${\;\;\;}$Now we would like to write down the formula for the covariantised 
product corresponding to a Drinfel'd double $D(A)$. It is slightly 
easier to do this  
using the explicit formulae obtained above for $Q$ and $Q^{-1}$ rather than 
the definition in terms of $\cdot $. We have then:
$$(a\otimes p)\underline{\cdot }(a'\otimes p')=Q^{-1}(Q(a\otimes p)
Q(a'\otimes p'))$$
for $a,a'\in A$ and $p,p'\in A^*$, and by a direct computation we obtain:
$$(a\otimes p)\underline{\cdot }(a'\otimes p')=\sum aa'_2\otimes 
(S^{-1}(a'_1)\rightharpoonup p\leftharpoonup a'_3)p'$$
and this is just the multiplication of the realization of the Drinfel'd 
double on $A\otimes A^*$, see \cite{maj1}, p. 290. Moreover, if 
$R=\sum (\varepsilon \otimes e_i)\otimes (e^i\otimes 1)$ is the $R$-matrix 
of $D(A)=A^*\otimes A$, then one can see that $(Q^{-1}\otimes Q^{-1})(R)=
\sum (e_i\otimes \varepsilon )\otimes (1\otimes e^i)$, which is the 
$R$-matrix of $D(A)=A\otimes A^*$. Finally, the comultiplication of 
$D(A)=A\otimes A^*$ is 
$$\Delta (a\otimes p)=\sum a_1\otimes p_2\otimes a_2\otimes p_1$$
and one can see that $Q:D(A)=A\otimes A^*\rightarrow D(A)=A^*\otimes A$ is a 
coalgebra map. \\
${\;\;\;}$In conclusion, for a Drinfel'd double, the map $Q$ in the 
definition of factorizability gives an isomorphism of quasitriangular Hopf 
algebras between the two realizations of the double, on $A^*\otimes A$ 
and $A\otimes A^*$. For the general meaning of the map $Q$ (i.e. for $any$ 
quasitriangular Hopf algebra) in terms of ``braided groups'', we refer to 
\cite{maj1}, p. 490.
\begin{re}{\em 
If we take the above realization of the double on $A\otimes A^*$ (denoted  
also by $D(A)$), then one can check that the multiplication of 
$D(A)^*_R$ is given by
$$(p\otimes a)\cdot (p'\otimes a')=\sum p'(a'_1\rightharpoonup p)\otimes 
a'_2a$$
so that $D(A)^*_R=$$\cal H$$(A^*)^{op}$. So, we can obtain $\cal H$$(A^*)$ 
from $D(A)$ by the same method used for obtaining $\cal H$$(A)$ (let us 
mention that in \cite{lu1} $\cal H$$(A^*)$ is obtained form $D(A)$ 
$without$ using the $R$-matrix of $D(A)$). }
\end{re}
${\;\;\;}$Now we shall speak about the pentagon equation, which, as the 
Yang-Baxter equation, appears in various contexts in mathematics and 
physics (see for instance \cite{bs}). If $A$ is a finite dimensional Hopf 
algebra, there exist two canonical procedures to construct (invertible) 
solutions for the pentagon equation: on the one hand, the map 
$w\in End (A\otimes A)$ given by $w(a\otimes b)=\sum a_1\otimes a_2b$ 
satisfies the pentagon equation $w_{12}w_{13}w_{23}=w_{23}w_{12}$, and 
its inverse is given by $w^{-1}(a\otimes b)=a_1\otimes S(a_2)b$ (see 
\cite{maj1}, p. 29). On the other hand, if we consider the element 
$W=\sum (1\otimes e^i)\otimes (e_i\otimes \varepsilon )$ in 
$\cal H$$(A)\otimes $$\cal H$$(A)$, where $\{e_i\}$ is a basis of $A$ and 
$\{e^i\}$ its dual basis in $A^*$, then $W$ is also a solution to the 
pentagon equation (see \cite{vdvk}, \cite{ka}). These two approaches are 
actually equivalent, because if we consider the algebra isomorphism 
$\lambda :A\#A^*\rightarrow End (A)$ given by $\lambda (a\otimes f)(b)=
\sum af(b_2)b_1$, for all $a,b\in A$ and $f\in A^*$ (see \cite{mon}, p. 162) 
then we have $(\lambda \otimes \lambda )(W)=w$. \\
${\;\;\;}$Since the Heisenberg double $\cal H$$(A)$ may be obtained from 
the Drinfel'd double $D(A)$, it is natural to see whether we can obtain the  
element $W$ from the $R$-matrix of $D(A)$. 
\begin{pr} If $A$ is a finite dimensional Hopf algebra, $R$ is the $R$-matrix 
of $D(A)$, $W$ is the canonical element of $\cal H$$(A)$ and 
$Q:D(A)^*\rightarrow D(A)$ is the map expressing the factorizability of 
$D(A)$, then we have 
$$W=(Q^{-1}\otimes Q^{-1})(R_{21})$$
where we identified $D(A)^*$ and $\cal H$$(A)$ as linear spaces. 
\end{pr}
{\bf Proof:} recall that the map $Q$ is given by $Q(a\otimes f)=\sum 
a_1\rightharpoonup f\leftharpoonup S^{-1}(a_3)\otimes a_2$, so we shall 
compute:\\[2mm]
$(Q\otimes Q)(W)=\sum Q(1\otimes e^i)\otimes Q(e_i\otimes \varepsilon )\\[2mm]
=\sum (e^i\otimes 1)\otimes ((e_i)_1\rightharpoonup \varepsilon  
\leftharpoonup S^{-1}((e_i)_3)\otimes (e_i)_2))\\[2mm]
=\sum (e^i\otimes 1)\otimes (\varepsilon \otimes e_i)=R_{21}$, q.e.d.\qed
${\;\;\;}$Let us note that although this proposition expresses $W$ in terms 
of the structure of $D(A)$, it does not give an $explanation$, in terms of 
$D(A)$, for $why$ is $W$ a solution to the pentagon equation. This would 
have been the case if we could have this proposition as a particular case of a 
more general statement, for instance, it was tempting to conjecture the 
following: ``If $(H, R)$ is a finite dimensional quasitriangular 
factorizable Hopf algebra, then the element 
$W=(Q^{-1}\otimes Q^{-1})(R_{21})$ is a solution to the pentagon equation 
with respect to the algebra structure of $H^*_R$''. Unfortunately, at least 
in this generality, this is false. A counterexample may be obtained  
as follows. Suppose $char (k)=0$ and $k$ contains a primitive root of 
unity of order 3, $\omega $ say. Take $G=\{1, a, a^2\}$ a group of order 3, 
and take $H=kG$, with the following $R$-matrix (see \cite{rad2}): 
$$R=\frac{1}{3}(1\otimes 1+1\otimes a+1\otimes a^2
+a\otimes 1+\omega ^2(a\otimes a)+
\omega (a\otimes a^2)+a^2\otimes 1+\omega (a^2\otimes a)+\omega ^2(a^2
\otimes a^2))$$
It is known from \cite{rad2} that $(H,R)$ is factorizable, a fact which may 
be proved also directly: if $\{e^0, e^1, e^2\}$ is the basis of $H^*$ 
dual to $\{1, a, a^2\}$, then one can check that the map $Q:H^*\rightarrow 
H$ is given by: 
$$Q(e^0)=\frac{1}{3}(1+a+a^2)$$
$$Q(e^1)=\frac{1}{3}(1+\omega a+\omega ^2a^2)$$
$$Q(e^2)=\frac{1}{3}(1+\omega ^2a+\omega a^2)$$
and that it is bijective, its inverse being the map $Q^{-1}:H\rightarrow H^*$,
$$Q^{-1}(1)=e^0+e^1+e^2$$
$$Q^{-1}(a)=e^0+\omega ^2e^1+\omega e^2$$
$$Q^{-1}(a^2)=e^0+\omega e^1+\omega ^2e^2$$
The multiplication of $H^*_R$ is given by:
$$e^0\cdot e^0=\frac{1}{3}(e^0+\omega ^2e^1+\omega ^2e^2)$$      
$$e^1\cdot e^1=\frac{1}{3}(e^1+\omega ^2e^0+\omega ^2e^2)$$
$$e^2\cdot e^2=\frac{1}{3}(e^2+\omega ^2e^0+\omega ^2e^1)$$
$$e^0\cdot e^1=e^1\cdot e^0=\frac{1}{3}(e^0+e^1+\omega e^2)$$
$$e^0\cdot e^2=e^2\cdot e^0=\frac{1}{3}(e^0+e^2+\omega e^1)$$
$$e^1\cdot e^2=e^2\cdot e^1=\frac{1}{3}(e^1+e^2+\omega e^0)$$
The element $W$ is given by 
$$W=(Q^{-1}\otimes Q^{-1})(R_{21})=$$
$$=e^0\otimes e^0+e^0\otimes e^1+
e^0\otimes e^2+e^1\otimes e^0+$$
$$+\omega (e^1\otimes e^1)+
\omega ^2(e^1\otimes e^2)+e^2\otimes e^0+\omega ^2(e^2\otimes e^1)+
\omega (e^2\otimes e^2)$$
and one can prove by a direct (but tedious) computation that $W$ is $not$ 
a solution to the pentagon equation. \\
${\;\;\;}$However, one can ask for what classes of Hopf algebras the 
conjecture is true. A useful result may be the following: 
\begin{pr} Let $(H, R)$ and $(H', R')$ be two quasitriangular Hopf algebras 
and $F:H\rightarrow H'$ a map of quasitriangular Hopf algebras. Then:\\
a) $F^*:H'^*_{R'}\rightarrow H^*_R$ is an algebra map\\
b) if $H$ and $H'$ are factorizable and $F$ is injective, then
$$(F^*\otimes F^*)((Q^{-1}\otimes Q^{-1})(R'_{21}))=
(Q^{-1}\otimes Q^{-1})(R_{21})$$
Consequently, if $(Q^{-1}\otimes Q^{-1})(R'_{21})$ is a solution to the 
pentagon equation, so is $(Q^{-1}\otimes Q^{-1})(R_{21})$.
\end{pr}
{\bf Proof:} a) has been noticed before, so we prove b). By applying 
$Q\otimes Q$, the relation we have to prove is equivalent to
$$\sum Q(F^*(Q'^{-1}(R'^1)))\otimes Q(F^*(Q'^{-1}(R'^2)))=R$$
which is equivalent to 
$$\sum Q'^{-1}(R'^1)(F(\alpha ^2)F(\beta ^1))\alpha ^1\beta ^2\otimes 
Q'^{-1}(R'^2)(F(\gamma ^2)F(\delta ^1))\gamma ^1\delta ^2=R$$
where $\alpha =\beta =\gamma =\delta =R$, which, by applying $F\otimes F$ 
(which is injective) and taking into account that $(F\otimes F)(R)=R'$, 
is equivalent to 
$$\sum Q'(Q'^{-1}(R'^1))\otimes Q'(Q'^{-1}(R'^2))=R'$$
which is obviously true.\qed
\begin{re}{\em 
Let $A$ be a finite dimensional Hopf algebra and define the map 
$\varphi :D(A)\rightarrow $$\cal H$$(A^*)\otimes $$\cal H$$(A)$ by 
$$\varphi (p\otimes a)=\sum (p_2\# a_1)\otimes (a_2\# p_1\leftharpoonup a_3)=
\sum (p_3\# a_1)\otimes (p_1\rightharpoonup a_2\# p_2)$$
for all $p\in A^*$ and $a\in A$. Then one can check, by a direct computation, 
that $\varphi $ is an algebra map (this is a ``coordinate-free'' version of a  
result in \cite{ka}) and that 
$$(id \otimes Q)\varphi =(id \otimes T\otimes id)\Delta _{D(A)}$$
where $Q:A\otimes A^*\rightarrow A^*\otimes A$, $Q(a\otimes p)=
\sum a_1\rightharpoonup p\leftharpoonup S^{-1}(a_3)\otimes a_2$ and 
$T:A\otimes A^*\rightarrow A\otimes A^*$, $T(a\otimes p)=\sum a_1\otimes 
a_2\rightharpoonup p$. \\
${\;\;\;}$Then, by identifying $\cal H$$(A^*)\simeq End(A^*)$ and 
$\cal H$$(A)\simeq End(A)$ via the algebra map $\lambda $ mentioned before 
and then by identifying $End (A^*)\otimes End (A)\simeq End (A^*\otimes A)$, 
and finally by composing with $\varphi $, we obtain an algebra map 
$$\alpha :D(A)\rightarrow End(A^*\otimes A)$$
$$\alpha (p\otimes a)(q\otimes b)=\sum p_2(a_1\rightharpoonup q)\otimes 
p_1\rightharpoonup (a_2b)$$
hence a $D(A)$-module structure on $A^*\otimes A$, which seems to be 
noncanonical.} 
\end{re}    
\section{A quantum groupoid with Heisenberg double as base}
${\;\;\;}$We recall first the following definitions from \cite{lu2}:
\begin{de} A bialgebroid consists of the following data: \\
1) an algebra $T$ called the total algebra\\
2) an algebra $B$ called the base algebra\\
3) the source map: an algebra homomorphism 
$$\alpha :B\rightarrow T$$
and the target map: an algebra anti-homomorphism
$$\beta :B\rightarrow T$$
such that $\alpha (b)\beta (b')=\beta (b')\alpha (b)$ for all $b, b'\in B$. \\
4) the coproduct: a $B-B$-bimodule map 
$$\Delta :T\rightarrow T\otimes _BT,\;\;h\mapsto \sum h_1\otimes h_2$$
with $\Delta (1)=1\otimes 1$, satisfying the coassociativity condition 
$$(\Delta \otimes _B id_T)\Delta =(id _T\otimes _B \Delta )\Delta :
T\rightarrow T\otimes _BT\otimes _BT$$
where $T$ is a $B-B$-bimodule with structure maps 
$$\lambda :B\otimes T\rightarrow T,\;\;b\otimes t\mapsto \alpha (b)t$$
$$\rho :T\otimes B\rightarrow T,\;\;t\otimes b\mapsto \beta (b)t$$
The coproduct $\Delta $ and the algebra structure of $T$ are required to be 
compatible in the sense that the kernel of the map 
$$\phi :T\otimes T\otimes T\rightarrow T\otimes _BT,\;\;t\otimes t'\otimes 
t''\mapsto (\Delta (t))(t'\otimes t'')$$
is a left ideal of $T\otimes T^{op}\otimes T^{op}$ (we have used the fact that 
$T\otimes T$ acts on $T\otimes _BT$ from the right by right multiplication).\\
5) the counit map: a $B-B$-bimodule map 
$$\varepsilon: T\rightarrow B$$
satisfying $\varepsilon (1)=1$ (it follows then that $\varepsilon \beta =
\varepsilon \alpha =id_B$) and  
$$\lambda (\varepsilon \otimes id_T)\Delta =\rho (id_T\otimes \varepsilon )
\Delta =id_T:T\rightarrow T$$
where $\lambda $ and $\rho $ are the maps defined above. \\
It is also required that $\varepsilon $ is compatible with the algebra 
structure of $T$ in the sense that the kernel of $\varepsilon $ is a left 
ideal of $T$. 
\end{de}  
\begin{re} {\em It was proved in \cite{lu2} that in the case 
$B=k$ the definition  
of a bialgebroid reduces to the one of a bialgebra over $k$.} 
\end{re}
\begin{re}{\em 
It was proved by P. Xu in \cite{xu} that the compatibility condition between 
$\Delta $ and the algebra structure of $T$ in the definition of a bialgebroid 
may be replaced by a (more manageable) equivalent one, namely the 
following two conditions have to be satisfied: 
$$\Delta (t)(\beta (b)\otimes 1-1\otimes \alpha (b))=0 \;\:in\;\;
T\otimes _BT$$
$$\Delta (tt')=\Delta (t)\Delta (t')\;\;in\;\;T\otimes _BT$$
for all $t, t'\in T$ and $b\in B$.} 
\end{re}
\begin{de} A Hopf algebroid (=quantum groupoid) is a bialgebroid $T$ over 
an algebra $B$ with structure maps $\alpha , \beta , \Delta \;and\; 
\varepsilon $ together with a linear map $\tau :T\rightarrow T$, called the 
antipode map, having the following properties: \\
1) $\tau $ is an algebra anti-isomorphism of $T$\\
2) $\tau \beta =\alpha $\\
3) $m_T(\tau \otimes id)\Delta =\beta \varepsilon \tau : T\rightarrow T$, 
where $m_T$ is the multiplication of $T$\\
4) there exists a linear map $\gamma :T\otimes _BT\rightarrow T\otimes T$ 
with the following properties:\\
4a) $\gamma $ is a section for the natural projection $p:T\otimes T
\rightarrow T\otimes _BT$\\
4b) the following identity holds:
$$m_T(id\otimes \tau )\gamma \Delta =\alpha \varepsilon :T\rightarrow T$$ 
\end{de}
\begin{re}{\em a) The existence of the section $\gamma $ is required because  
the map $m_T(id\otimes \tau )\Delta $ is $not$ well-defined.\\
b) in general $\varepsilon \tau \neq \varepsilon $, see \cite{lu2}.}
\end{re}
${\;\;\;}$Recall also from \cite{lu2} the following results:
\begin{te} (\cite{lu2}) If $A$ is a finite dimensional Hopf algebra and $V$  
is a quantum commutative left $D(A)$-module algebra, then $V\#A$ becomes a  
Hopf algebroid over $A$, the left action of $A$ on $V$ being the one 
induced from the action of $D(A)$ on $V$ by the inclusion $A\subseteq D(A)$.
\end{te}
\begin{pr} (\cite{lu2}) If $A$ is a finite dimensional Hopf algebra then 
$A^*$ is a quantum commutative left $D(A)$-module algebra, so the Heisenberg 
double $A^*\#A$ of $A^*$ is a Hopf algebroid over $A^*$. 
\end{pr}
${\;\;\;}$Now, we want to obtain an example of a Hopf algebroid having the 
Heisenberg double of a finite dimensional Hopf algebra $A$ as base. We 
proceed as follows: first we give a generalization of the above theorem and 
then we use the quantum commutativity of the Heisenberg double of $A$ as a 
$D(A)-D(A)^{cop}$-bimodule algebra. \\
${\;\;\;}$Let us recall first some notation and results from \cite{rad1}. 
If $(H,R)$ is a quasitriangular Hopf algebra, define the subspaces 
$L=R_{(l)}$ and $D=R_{(r)}$ by $R_{(l)}=\{(id\otimes p)(R)/p\in H^*\}$ and  
$R_{(r)}=\{(p\otimes id)(R)/p\in H^*\}$. If we write $R=\sum _{i=1}^mu_i 
\otimes v_i\in H\otimes H$ where $m$ is as small as possible, then 
$\{u_1,...,u_m\}$ is a basis for $L$ and $\{v_1,...,v_m\}$ is a basis 
for $D$, in particular $dim L=dim D$ and this common dimension is 
called the rank of $R$ and is denoted by $rank (R)$. Moreover, $L$ and $D$ 
are (finite-dimensional) Hopf subalgebras of $H$ and the map $f:L^{*cop}
\rightarrow D$ defined by $f(p)=(p\otimes id)(R)$ for $p\in L^*$ is a 
Hopf algebra isomorphism. \\
${\;\;\;}$Now we can generalize Lu's theorem as follows:    
\begin{te} Let $(H,R)$ be a quasitriangular Hopf algebra and denote as above 
$L=R_{(l)}$, $D=R_{(r)}$. Let $V$ be a quantum commutative left $H$-module 
algebra with action denoted by $H\otimes V\rightarrow V$, $h\otimes v\mapsto 
h\cdot v$. Then $V\#L$ is a Hopf algebroid over $V$, the left action of $L$ 
on $V$ being the one induced from the action of $H$ on $V$ by the inclusion 
$L\subseteq H$. 
\end{te}
{\bf Proof:} since (with one exception) the proof follows closely the 
proof in \cite{lu2}, we  
shall just write down the structure maps and we shall give only those parts 
of the proof which, in \cite{lu2}, depend on the explicit structure of the 
Drinfel'd double. \\
The map $\alpha :V\rightarrow V\#L$ is given by $\alpha (v)=v\#1$.\\
Then, we know that the map $f:L^{*cop}\rightarrow D$, $f(p)=\sum p(R^1)R^2$  
is an isomorphism of Hopf algebras, so $V$ becomes a left $L^{*cop}$-module 
algebra via $f$, in particular it becomes a right $L$-comodule; denote by 
$\beta :V\rightarrow V\otimes L$ the structure map. If we denote by 
$\beta (v)=\sum v_0\otimes v_1$, for $v\in V$, we know that this is 
equivalent to $p\cdot v=\sum p(v_1)v_0$ for all $p\in L^*$, where $p\cdot v=
f(p)\cdot v=\sum p(R^1)R^2\cdot v$. So, $\beta $ is given by the formula
$$\beta (v)=\sum R^2\cdot v\otimes R^1$$
Then, exactly as in \cite{lu2}, using the properties of the $R$-matrix $R$ 
and the quantum commutativity of $V$, one can prove that $\beta $ is an 
algebra anti-homomorphism and that $\alpha (v)\beta (w)=\beta (w)\alpha (v)$ 
for all $v, w\in V$. \\
Let us note that, for $v, w\in V$ and $l\in L$, we have in $V\#L$:\\[2mm]
$\beta (v)(w\#l)=\sum (R^2\cdot v\#R^1)(w\#l)\\[2mm]
=\sum (R^2\cdot v)((R^1)_1\cdot w)\#(R^1)_2l\\[2mm]
=\sum (R^2r^2\cdot v)(R^1\cdot w)\#r^1l\\[2mm]
=\sum w(r^2\cdot v)\#r^1l$\\[2mm]
(we have used the quantum commutativity of $V$ and the fact that 
$(\Delta \otimes id)(R)=\sum R^1\otimes r^1\otimes R^2r^2$)\\[2mm]
The $V-V$-bimodule structure of $V\#L$ is given by 
$$v\cdot (w\#l)=\alpha (v)(w\#l)=(v\#1)(w\#l)=vw\#l$$
$$(w\#l)\cdot v=\beta (v)(w\#l)=\sum w(r^2\cdot v)\#r^1l$$
The coproduct is given by
$$\Delta :V\#L\rightarrow (V\#L)\otimes _V(V\#L)$$
$$\Delta (v\#l)=\sum (v\#l_1)\otimes _V(1\#l_2)$$
Obviously $\Delta (1\#1)=(1\#1)\otimes _V(1\#1)$ and $\Delta $ is a morphism  
of left $V$-modules. Let us prove that it is also a morphism of right 
$V$-modules. We calculate: \\[2mm]
$\Delta ((w\#l)\cdot v)=\Delta (\sum w(r^2\cdot v)\#r^1l)\\[2mm]
=\sum (w(r^2\cdot v)\#(r^1)_1l_1)\otimes _V(1\#(r^1)_2l_2)\\[2mm]
=\sum (w(R^2r^2\cdot v)\#R^1l_1)\otimes _V(1\#r^1l_2)$\\[3mm]
$\Delta (w\#l)\cdot v=\sum (w\#l_1)\otimes _V(1\#l_2)\cdot v\\[2mm]
=\sum (w\#l_1)\otimes _V(r^2\cdot v\#r^1l_2)\\[2mm]
=\sum (w\#l_1)\otimes _V(r^2\cdot v)\cdot (1\#r^1l_2)\\[2mm]
=\sum (w\#l_1)\cdot (r^2\cdot v)\otimes _V(1\#r^1l_2)$\\[2mm]
(since the tensor product is over $V$)\\[2mm]
$=\sum (w(R^2r^2\cdot v)\#R^1l_1)\otimes _V(1\#r^1l_2)$, q.e.d.\\
The coassociativity of $\Delta $ follows immediately from the one of 
$\Delta _L$. \\
The compatibility of $\Delta $ with the product of $V\#L$ may be proved as 
in \cite{lu2}, once we prove the following identity which in \cite{lu2} is 
proved using the structure of the Drinfel'd double:
$$\sum \beta (l_2\cdot v)(1\#l_1)=(1\#l)\beta (v)$$
for all $v\in V,l\in L$. To prove this we compute:\\[2mm]
$\sum \beta (l_2\cdot v)(1\#l_1)=\sum (R^2l_2\cdot v\#R^1)(1\#l_1)=\sum 
R^2l_2\cdot v\#R^1l_1$\\[3mm]
$(1\#l)\beta (v)=\sum (1\#l)(R^2\cdot v\#R^1)\\[2mm]
=\sum l_1R^2\cdot v\#l_2R^1\\[2mm]
=\sum R^2l_2\cdot v\#R^1l_1$\\[2mm]
using the identity $\Delta ^{cop}(l)R=R\Delta (l)$.\\
The section $\gamma :(V\#L)\otimes _V(V\#L)\rightarrow (V\#L)\otimes (V\#L)$ 
which is needed in this part of the proof (and also at the end) is given as 
in \cite{lu2} by  
$$\gamma ((v\#l)\otimes _V(v'\#l'))=\beta (v')(v\#l)\otimes (1\#l')$$
${\;\;\;}$Since the proof of this part in \cite{lu2} is rather tricky, 
let us give a proof using Xu's equivalent condition. We have to prove first 
that 
$$\Delta (v\#l)(\beta (w)\otimes 1-1\otimes \alpha (w))=0\;\;in\;\;
(V\#L)\otimes _V(V\#L)$$
for all $v, w\in V$ and $l\in L$. We calculate:\\[2mm]
$\Delta (v\#l)(\beta (w)\otimes 1-1\otimes \alpha (w))=\\[2mm]
=\sum ((v\#l_1)\otimes _V(1\#l_2))(\sum (R^2\cdot w\#R^1)\otimes (1\#1) 
-(1\#1)\otimes (w\#1))\\[2mm]
=\sum (v\#l_1)(R^2\cdot w\#R^1)\otimes _V(1\#l_2)-\sum (v\#l_1)\otimes _V
(1\#l_2)(w\#1)\\[2mm]
=\sum (v(l_1\cdot (R^2\cdot w))\#l_2R^1)\otimes _V(1\#l_3)-
\sum (v\#l_1)\otimes _V(l_2\cdot w\#l_3)\\[2mm]
=\sum (v(l_1\cdot (R^2\cdot w))\#l_2R^1)\otimes _V(1\#l_3)-\sum (v\#l_1)
\otimes _V(l_2\cdot w)\cdot (1\#l_3)\\[2mm]
=\sum (v(l_1\cdot (R^2\cdot w))\#l_2R^1)\otimes _V(1\#l_3)-
\sum (v\#l_1)\cdot (l_2\cdot w)\otimes _V(1\#l_3)$\\[2mm]
(since the tensor product is over $V$)\\[2mm]
$=\sum (v(l_1\cdot (R^2\cdot w))\#l_2R^1)\otimes _V(1\#l_3)-
\sum (v(R^2\cdot (l_2\cdot w))\#R^1l_1)\otimes _V(1\#l_3)$\\[2mm]
(we used the formula for the right module structure of $V\#L$)\\[2mm]
$=\sum (v(l_1R^2\cdot w)\#l_2R^1)\otimes _V(1\#l_3)-
\sum (v(R^2l_2\cdot w)\#R^1l_1)\otimes _V(1\#l_3)$\\[2mm]
and this is zero because $\Delta ^{cop}(l)R=R\Delta (l)$. \\
Next we have to prove that 
$$\Delta ((v\#l)(v'\#l'))=\Delta (v\#l)\Delta (v'\#l')$$
in $(V\#L)\otimes _V(V\#L)$ for all $v, v'\in V$ and $l,l'\in L$, and 
this is trivial. \\
${\;\;\;}$The counit is given by 
$$\varepsilon :V\#L\rightarrow V$$
$$\varepsilon (v\#l)=\varepsilon (l)v$$
and it is easy to see that it satisfies all the required properties. \\
So, we have obtained so far a bialgebroid structure on $V\#L$ over $V$. We 
shall construct its antipode. Let $u=\sum S(R^2)R^1$ be the canonical 
Drinfel'd element of $H$ and define $d_0=S(u)^{-1}$ (this element is denoted 
by $u_2$ in \cite{dr1}). Then we know from \cite{dr1} that $d_0$ has the 
following properties:\\[2mm]
$d_0=\sum S^2(R^1)R^2$\\[2mm]
$d_0^{-1}=\sum S^{-1}(R^1)R^2=\sum R^1S(R^2)$\\[2mm]
$S^2(h)=d_0hd_0^{-1}$ for all $h\in H$\\[2mm]
$\Delta (d_0)=(R_{21}R)(d_0\otimes d_0)=(d_0\otimes d_0)(R_{21}R)$\\[2mm]
As in \cite{lu2} one can prove that $d_0$ acts as an algebra isomorphism on  
$V$.\\
Define now the antipode of $V\#L$, by
$$\tau :V\#L\rightarrow V\#L$$
$$\tau (v\#l)=(1\#S(l))\beta (d_0\cdot v)$$
We shall prove that $\tau \beta =\alpha $. Let $v\in V$; we calculate: \\[2mm]
$\tau (\beta (v))=\sum \tau (R^2\cdot v\#R^1)=\sum (1\#S(R^1))\beta 
(d_0R^2\cdot v)\\[2mm]
=\sum (1\#S(R^1))(r^2d_0R^2\cdot v\#r^1)\\[2mm]
=\sum S((R^1)_2)r^2d_0R^2\cdot v\#S((R^1)_1)r^1\\[2mm]
=\sum S(\rho ^1)r^2d_0R^2\rho ^2\cdot v\#S(R^1)r^1$\\[2mm]
(where $\rho =R$ and we used the relation $(\Delta \otimes id)(R)=\sum R^1
\otimes \rho ^1\otimes R^2\rho ^2$) \\[2mm]
$=\sum S(\rho ^1)r^2S^2(R^2)d_0\rho ^2\cdot v\#S(R^1)r^1$\\[2mm]
(using $S^2(x)=d_0xd_0^{-1}$)\\[2mm]
$=\sum S(\rho ^1)r^2S(R^2)d_0\rho ^2\cdot v\#R^1r^1$\\[2mm]
(since $(S\otimes S)(R)=R$, see \cite{dr1})\\[2mm]
$=\sum S(\rho ^1)d_0\rho ^2\cdot v\#1$\\[2mm]
(since $R^{-1}=\sum R^1\otimes S^{-1}(R^2)$, see \cite{dr1})\\[2mm]
$=\sum S(\rho ^1)S^2(\rho ^2)d_0\cdot v\#1$\\[2mm]
(again using $S^2(x)=d_0xd_0^{-1}$)\\[2mm]
$=\sum \rho ^1S(\rho ^2)d_0\cdot v\#1$\\[2mm]
(again since $(S\otimes S)(\rho )=\rho $)\\[2mm]
$=v\#1$\\[2mm]
(since $d_0^{-1}=\sum \rho ^1S(\rho ^2)$)\\[2mm]
$=\alpha (v)$, q.e.d.\\
The inverse of $\tau $ is given by 
$$\tau ^{-1}:V\#L\rightarrow V\#L$$
$$\tau ^{-1}(v\#l)=(1\#S^{-1}(l))\beta (v)$$
One can prove by a direct computation that  
$\tau ^{-1}$ is an algebra anti-homomorphism, then using this one can prove 
that $\tau ^{-1}\tau =id$ and finally that $\tau \tau ^{-1}=id$ by a 
computation similar to the proof of $\tau \beta =\alpha $.  
Hence $\tau $ is an algebra  
anti-isomorphism. The other axioms for $\tau $ are proved exactly as in  
\cite{lu2}. \qed 
\begin{re} {\em 
a) the proof shows that the statement in the theorem remains valid if we 
replace $L$ by any Hopf subalgebra of $H$ containing $L$.\\ 
b) if $A$ is a finite dimensional Hopf algebra and $(D(A), R)$ its 
Drinfel'd double, it is easy to see that $R_{(l)}=\{\varepsilon 
\otimes a\;/\;a\in A\}$, which is isomorphic as Hopf algebras to $A$, so  
indeed the above result generalizes Lu's theorem.\\
c) one can check that, in the case $H=D(A)$ and $V=A^*$, the formula  
$\tau ^{-1}(v\#l)=(1\#S^{-1}(l))\beta (v)$ for the inverse of the antipode 
of $\cal H$$(A^*)$ may be written as
$$\tau ^{-1}(p\#a)=\sum (S^{-1}(a)\rightharpoonup (e^ip))e^j\#S^{-1}(e_j)e_i$$
for all $a\in A$ and $p\in A^*$, where $\{e_i\}$ is a basis in $A$ and 
$\{e^i\}$ its dual basis in $A^*$ (the formula for $\tau $ was given in 
\cite{lu2}).\\ 
d) the generalization of Lu's theorem in \cite{bm} (independently obtained, 
as we mentioned before) is more general than ours: it concerns quantum 
commutative algebras in categories of Yetter-Drinfel'd modules.} 
\end{re}   
${\;\;\;}$Now, let $A$ be a finite dimensional Hopf algebra. As we have seen 
before, $\cal H$$(A)$ is a quantum commutative left $D(A)\otimes 
D(A)^{op\;cop}$-module algebra, so the above theorem may be applied. 
Denote by $\cal R$$=\sum (R^1\otimes U^2)\otimes (R^2\otimes U^1)$ the 
$R$-matrix of $D(A)\otimes D(A)^{op\;cop}$, where $R$ is the $R$-matrix of 
$D(A)$ and $U=R^{-1}$. Since $R^{-1}=\sum S(R^1)\otimes R^2$ (see \cite{dr1}) 
and the antipode of $D(A)$ has the property that $S_{D(A)}(\varepsilon 
\otimes a)=\varepsilon \otimes S(a)$ for all $a\in A$, it follows that  
$\cal R$ is given by $\cal R$$=\sum (\varepsilon \otimes e_i\otimes e^j
\otimes 1)\otimes (e^i\otimes 1\otimes \varepsilon \otimes S(e_j))$, where 
$\{e_i\}$ is a basis in $A$ and $\{e^i\}$ its dual basis in $A^*$. Then it 
is easy to see that the set $\{\varepsilon \otimes e_i\otimes e^j\otimes 1\}$ 
is a $k$-linear basis in $\cal R$$_{(l)}$, so $\cal R$$_{(l)}$ may be 
identified as linear spaces with $A\otimes A^*$, and from the Hopf algebra 
structure of $D(A)\otimes D(A)^{op\;cop}$ we can see that actually  
$\cal R$$_{(l)}$ may be identified with $A\otimes A^{*op}$ as Hopf 
subalgebras of $D(A)\otimes D(A)^{op\;cop}$. Some calculations with the  
explicit formulae given before for the left and right regular actions of 
$D(A)$ on $\cal H$$(A)$ show that the left $A\otimes A^{*op}$-module 
structure of $\cal H$$(A)$ is given by:
$$(a\otimes p)\cdot (b\#q)=(b\leftharpoonup p)\#(a\rightharpoonup q)=
\sum p(b_1)q_2(a)(b_2\#q_1)$$    
for all $a, b\in A$ and $p, q\in A^*$. In conclusion, we obtain:
\begin{co} $\cal H$$(A)\# (A\otimes A^{*op})$ is a Hopf algebroid over 
$\cal H$$(A)$.
\end{co}  
\section{Something like a vertex group}
${\;\;\;}$We start with the following definition, which is a variation of 
the one in \cite{b1}. 
\begin{de} Let $H$ be a cocommutative Hopf algebra. A vertex group over 
$H$ is an $H$-bimodule algebra $K$ (with actions denoted by $h\otimes k
\mapsto h\cdot k$ and $k\otimes h\mapsto k\cdot h$ for all $h\in H$ 
and $k\in K$) together with a map $\alpha :H^*\rightarrow K$ which is a 
morphism of $H$-bimodule algebras ($H^*$ is an $H$-bimodule algebra via 
the left and right regular actions) and an algebra anti-isomorphism 
$\tau :K\rightarrow K$ satisfying the conditions:
$$\tau (\alpha (p))=\alpha (S(p))$$
$$\tau (h\cdot k)=\tau (k)\cdot S(h)$$
$$\tau (k\cdot h)=S(h)\cdot \tau (k)$$
for all $h\in H, p\in H^*, k\in K$, where we denoted by $S$ the antipode 
of $H$ and also the antipode of $H^*$. The map $\tau $ will be called 
the antipode of $K$. If $K$ is commutative and $\tau ^2=id$ we recover 
Definition 3.2 in \cite{b1}.
\end{de} 
\begin{re}{\em Of course, $H^*$ is a vertex group over $H$ with $\alpha =id$  
and $\tau =S$.}
\end{re}
\begin{pr} If $A$ is a finite dimensional cocommutative Hopf algebra, 
then $\cal H$$(A)$ is a vertex group over $A$. 
\end{pr}
{\bf Proof:} define $\alpha :A^*\rightarrow $$\cal H$$(A)=A\#A^*$, 
$\alpha (p)=1\#p$, which is of course an algebra map. We have seen that 
$\cal H$$(A)$ is a $D(A)-D(A)^{cop}$-bimodule algebra; by restricting  
the actions of $D(A)$ on $\cal H$$(A)$ to actions of $A$ and since  
$A$ is cocommutative, we obtain that $\cal H$$(A)$ is an $A$-bimodule  
algebra. \\
Take $\{e_i\}$ a basis in $A$ and $\{e^i\}$ its dual basis in $A^*$ and define 
$\tau :A\#A^*\rightarrow A\#A^*$ by
$$\tau (a\#p)=\sum (S(p)\rightharpoonup (e_ia))e_j\#S(e^j)e^i$$
This is just the map $\tau ^{-1}$ in the previous section 
written for $\cal H$$(A)$   
instead of $\cal H$$(A^*)$ (we have used also the fact that $S^2=id$); 
we already know that this is an algebra anti-isomorphism.\\
${\;\;\;}$Let us check the other axioms of a vertex group. In order to 
simplify the calculations we shall use the version of Sweedler's sigma 
notation introduced in \cite{mas} for dealing with cocommutative Hopf 
algebras, namely for $a\in A$ we shall denote $\Delta (a)=\sum a\otimes a$, 
$(id\otimes \Delta )\Delta (a)=\sum a\otimes a\otimes a$ etc. With this 
notation, the antipode axiom for $A$ is written as $\sum S(a)a=\sum aS(a)=
\varepsilon (a)1$ for all $a\in A$. \\
${\;\;\;}$Let us prove that $\alpha $ is a bimodule map, namely 
$$\alpha (a\rightharpoonup p\leftharpoonup b)=a\rightharpoonup \alpha (p)
\leftharpoonup b$$
for all $a,b\in A$ and $p\in A^*$, that is 
$$1\#a\rightharpoonup p\leftharpoonup b=(\varepsilon \otimes a)
\rightharpoonup (1\#p)\leftharpoonup (\varepsilon \otimes b)$$
where $\varepsilon \otimes a$ and $\varepsilon \otimes b$ are considered 
as elements in the Drinfel'd double of $A$. By evaluating in an element 
$\varphi \otimes x\in A^*\otimes A$, the lhs becomes $\varphi (1)p(bxa)$ 
and the rhs may be calculated as:\\[2mm]
$(1\#p)((\varepsilon \otimes b)(\varphi \otimes x)(\varepsilon 
\otimes a))\\[2mm]
=(1\#p)((\varepsilon \otimes b)(\varphi \otimes xa))\\[2mm]
=(1\#p)(b\rightharpoonup \varphi \leftharpoonup S(b)\otimes bxa)\\[2mm]
=\varphi (S(b)b)p(bxa)=\varphi (1)p(bxa)$, q.e.d.\\[2mm]
It is very easy to see that $\tau (1\#p)=1\#S(p)$, so let us check the last 
two axioms for $\tau $. From the general formulae for the actions of $D(A)$ 
on $\cal H$$(A)$ given in a previous section we can write down the actions 
of $A$ on $\cal H$$(A)$: 
$$a\rightharpoonup (b\#p)=\sum p_2(a)(b\#p_1)$$
$$(b\#p)\leftharpoonup a=\sum p_1(a)(S(a)ba\#p_2)$$
for all $a,b\in A$ and $p\in A^*$. We calculate: \\[2mm]
$\tau (a\rightharpoonup (b\#p))=\sum p_2(a)\tau (b\#p_1)\\[2mm]
=\sum p_2(a)(S(p_1)\rightharpoonup (e_ib))e_j\#S(e^j)e^i$\\[2mm]
which evaluated in an element $\varphi \otimes x\in A^*\otimes A$ gives \\[2mm]
$\;\;\;\sum p(S(b)S(x)a)\varphi (xbS(x))$\\[2mm]
On the other hand, we have:\\[2mm]
$\tau (b\#p)\leftharpoonup S(a)=\\[2mm]
=\sum ((S(p)\rightharpoonup (e_ib))e_j\#S(e^j)e^i)\leftharpoonup S(a)\\[2mm]
=\sum S((e^j)_2)(S(a))(e^i)_1(S(a))(a(S(p)\rightharpoonup (e_ib))e_jS(a)\# 
S((e^j)_1)(e^i)_2)$\\[2mm]
which evaluated in $\varphi \otimes x\in A^*\otimes A$ gives:\\[2mm]
$\sum (e^j)_2(a)(e^i)_1(S(a))p(S(b)S(e_i))\varphi (ae_ibe_jS(a))
(e^j)_1(S(x))(e^i)_2(x)\\[2mm]
=\sum e^i(S(a)x)e^j(S(x)a)p(S(b)S(e_i))\varphi (ae_ibe_jS(a))\\[2mm]
=\sum p(S(b)S(x)a)\varphi (aS(a)xbS(x)aS(a))\\[2mm]
=\sum p(S(b)S(x)a)\varphi (xbS(x))$\\[2mm]
so we have obtained $\tau (a\rightharpoonup (b\#p))=\tau (b\#p)\leftharpoonup 
S(a)$. \\
${\;\;\;}$Now we prove the last relation. \\[2mm]
$\tau ((b\#p)\leftharpoonup a)=\sum p_1(a)\tau (S(a)ba\#p_2)\\[2mm]
=\sum p_1(a)((S(p_2)\rightharpoonup (e_iS(a)ba))e_j\#S(e^j)e^i)\\[2mm]
=\sum p_1(a)p_2(S(a)S(b)aS(e_i))(e_iS(a)bae_j\#S(e^j)e^i)\\[2mm]
=\sum p(S(b)aS(e_i))(e_iS(a)bae_j\#S(e^j)e^i)$\\[2mm]
which evaluated in $\varphi \otimes x\in A^*\otimes A$ gives\\[2mm] 
$\sum p(S(b)aS(x))\varphi (xS(a)baS(x))$\\[3mm]
$S(a)\rightharpoonup \tau (b\#p)=\sum p(S(b)S(e_i))S(a)\rightharpoonup 
(e_ibe_j\#S(e^j)e^i)\\[2mm]
=\sum p(S(b)S(e_i))S((e^j)_1)(S(a))(e^i)_2(S(a))(e_ibe_j\#S((e^j)_2)
(e^i)_1)$\\[2mm]
which evaluated in $\varphi \otimes x\in A^*\otimes A$ gives:\\[2mm]
$\sum p(S(b)S(e_i))e^i(xS(a))e^j(aS(x))\varphi (e_ibe_j)\\[2mm]
=\sum p(S(b)aS(x))\varphi (xS(a)baS(x))$\\[2mm]
so that $\tau ((b\#p)\leftharpoonup a)=S(a)\rightharpoonup \tau (b\#p)$, q.e.d.
 \qed

\section{The Heisenberg double for quasi-Hopf algebras}
${\;\;\;}$In this section we return to quasi-Hopf algebras and suggest a 
possible definition for Heisenberg doubles of them. \\
${\;\;\;}$If $A$ is a finite dimensional quasi-Hopf algebra, it was proved in 
\cite{maj2} that there exists a unique (up to isomorphism) quasitriangular  
quasi-Hopf algebra $(D(A), R)$ with the property that the category $D(A)$-mod 
is braided equivalent to the centre of the tensor category $A$-mod. This 
$D(A)$ is called the quantum double of $A$ (it generalizes the Drinfel'd 
double for Hopf algebras). Some concrete realizations of $D(A)$ (on 
$A^*\otimes A$ and $A\otimes A^*$) have been given in \cite{hn1}, 
\cite{hn2}.\\
${\;\;\;}$Now we can propose the following
\begin{de} The Heisenberg double of $A$ is $D(A)^*_R$, which is, as we know, 
a quantum commutative $D(A)-D(A)^{cop}$-bimodule algebra. It is denoted by  
$\cal H$$(A)$. 
\end{de}
${\;\;\;}$From a previous discussion and the uniqueness of $D(A)$, it 
follows that $\cal H$$(A)$ is ``unique up to isomorphism''. \\
${\;\;\;}$Because of the complicated structure of $D(A)$, it is quite 
difficult to write down explicitly the formula for the multiplication in 
$\cal H$$(A)$, and we shall not do this here. Instead, we shall do this in a 
particular case, related to the Dijkgraaf-Pasquier-Roche quasi-Hopf algebras 
$D^{\omega }(G)$ introduced in \cite{dpr}. 
We shall work in the slightly more general situation of \cite{bp}  
and we begin by recalling the relevant parts of the construction in 
\cite{bp}.    \\
${\;\;\;}$Let $H$ be a finite dimensional cocommutative
Hopf algebra (for which we shall use the variation of $\Sigma $-notation as 
in the previous section) with antipode $S$, and 
$\omega:H\otimes H\otimes H\rightarrow k$ a normalized
3-cocycle in the Sweedler cohomology \cite{sw1},  
that is $\omega $ is $k$-linear, convolution 
invertible and satisfies the conditions:
$$\sum \omega (x, y, zt)
\omega (xy, z, t)=\sum \omega (y, z, t)
\omega (x, yz, t)\omega (x, y, z)$$
$$\omega (1, x, y)=\omega (x, 1, y)=\omega (x, y, 1)=
\varepsilon (x)\varepsilon (y)$$ 
for all $x, y, z, t\in H$. Introduce also the following notation: 
$g\triangleleft x=\sum S(x)gx$ for all $g, x\in H$. Then, on the $k$-linear 
space $H^*\otimes H$ may be constructed a quasitriangular quasi-Hopf algebra, 
denoted by $D^{\omega }(H)$, for which the multiplication is given by
$$(p\otimes h)(p'\otimes h')=\sum p(h\rightharpoonup p'\leftharpoonup S(h))
\sigma (h, h')\otimes hh'$$
for all $p, p'\in H^*$ and $h, h'\in H$, where $\sigma :H\otimes H\rightarrow 
H^*$ is given by
$$\sigma (x, y)(g)=\theta (g; x, y)$$
for all $x, y, g\in H$, where $\theta :H\otimes H\otimes H\rightarrow k$,  
$$\theta (g; x, y)=\sum \omega (g, x, y)
\omega (x, y, g\triangleleft (xy))
\omega ^{-1}(x, g\triangleleft x, y)$$
for $g, x, y\in H$, where $\omega ^{-1}$ is the convolution inverse 
of $\omega $.\\
${\;\;\;}$Define the map $\gamma:H\otimes H\otimes H\rightarrow k$ by
$$\gamma (g, h; x)=\sum \omega (g, h, x)
\omega (x, g\triangleleft x, h\triangleleft x)\omega ^{-1}(g, x,
h\triangleleft x)$$
for $x, g, h\in H$, and the map $\nu: H\rightarrow (H\otimes H)^*$,
by $\nu (h)(x\otimes y)=\gamma (x, y; h)$.
Identifying $(H\otimes H)^*$ with $H^*\otimes H^*$, we shall write,
for any $h\in H$, $\nu (h)=\sum \nu _1(h)\otimes \nu _2(h)\in
H^*\otimes H^*$, and this relation is equivalent to
$\nu (h)(x\otimes y)=\sum \nu _1(h)(x)\nu _2(h)(y)$
for all $x, y\in H$. Then the comultiplication of $D^{\omega }(H)$ is given by 
$$\Delta : D^{\omega }(H)\rightarrow D^{\omega}(H)\otimes
D^{\omega }(H)$$
$$\Delta (p\otimes h)=
\sum (\nu _1(h)p_1\otimes h)\otimes (\nu _2(h)p_2\otimes h)$$
The associator of $D^{\omega }(H)$ is $\omega ^{-1}$ (regarded in 
$(H^*)^{\otimes 3}\subseteq D^{\omega }(H)^{\otimes 3}$) and the $R$-matrix is 
$$R=\sum (e^i\otimes 1)\otimes (\varepsilon \otimes e_i)$$
where $\{e_i\}$ is a basis in $H$ and $\{e^i\}$ its dual basis in $H^*$. 
Recall also from \cite{bp} the following identity: 
$$\sum \gamma (x, y\triangleleft x; h)
\theta (y; x, h)=\sum \theta (y; h, x\triangleleft h)
\gamma (y, x; h)$$
${\;\;\;}$It was proved in \cite{pvo} that $D^{\omega }(H)$ is (isomorphic 
to) the quantum double of a certain quasi-Hopf algebra denoted by 
$H^*_{\omega }$, which is just $H^*$ with its usual Hopf algebra structure 
but with a nontrivial associator, namely $\omega ^{-1}$. \\
${\;\;\;}$We can recover the DPR quasi-Hopf algebra $D^{\omega }(G)$ in 
\cite{dpr} by taking $H$ to be the group algebra of a finite group $G$.\\
${\;\;\;}$Let us note that if $\omega $ is trivial, then $D^{\omega }(H)$ 
is a Hopf algebra, isomorphic to the Drinfel'd double of $H^*$; but since 
we worked with the realization of $D^{\omega }(H)$ on $H^*\otimes H$ rather 
than $H\otimes H^*$, we can see from the explicit formulae given above that 
$(D^{triv}(H), R)=(D(H)^{cop}, R_{21})$, where $D(H)$ is the usual 
Drinfel'd double of $H$ and $R$ in the rhs is the $R$-matrix of $D(H)$.\\
${\;\;\;}$We can compute now the multiplication in $D^{\omega }(H)^*_R$ 
(which is, by our definition, the Heisenberg double of $H^*_{\omega }$). Take  
$x\otimes \varphi ,x'\otimes \varphi '\in D^{\omega }(H)^*_R$ and evaluate 
their product against an element $p\otimes h\in D^{\omega }(H)$:\\[2mm]
$((x\otimes \varphi )\cdot (x'\otimes \varphi '))(p\otimes h)=\\[2mm]
=\sum (x'\otimes \varphi ')((p\otimes h)_1R^2)(x\otimes \varphi )((p\otimes h)
_2R^1)\\[2mm]
=\sum (x'\otimes \varphi ')((\nu _1(h)p_1\otimes h)(\varepsilon \otimes e_i))
(x\otimes \varphi )((\nu _2(h)p_2\otimes h)(e^i\otimes 1))\\[2mm]
=\sum (x'\otimes \varphi ')(\nu _1(h)p_1\sigma (h, e_i)\otimes he_i)
(x\otimes \varphi )(\nu _2(h)p_2(h\rightharpoonup e^i\leftharpoonup S(h))
\otimes h)\\[2mm]
=\sum \nu _1(h)(x')p_1(x')\sigma (h, e_i)(x')\varphi '(he_i)\nu _2(h)(x)
p_2(x)e^i(S(h)xh)\varphi (h)\\[2mm]
=\sum p(x'x)\varphi '(xh)\varphi (h)\gamma (x', x; h)\theta (x'; h, 
x\triangleleft h)\\[2mm]
=\sum p(x'x)\varphi '(xh)\varphi (h)\gamma (x, x'\triangleleft x; h)
\theta (x'; x, h)$\\[2mm]
So, we have obtained the following formula:
$$(x\otimes \varphi )\cdot (x'\otimes \varphi ')=\sum x'(\varphi '_1
\rightharpoonup x)\otimes \varphi '_2\varphi \gamma (x, x'\triangleleft x; 
\cdot)\theta (x'; x, \cdot )$$
Let us note that if $\omega $ is trivial this formula becomes  
$$(x\otimes \varphi )\cdot (x'\otimes \varphi ')=\sum x'(\varphi '_1
\rightharpoonup x)\otimes \varphi '_2\varphi $$
which is the multiplication of $\cal H$$(H)^{op}$ (this is not a surprize, 
since $(D^{triv}(H), R)=(D(H)^{cop}, R_{21})$), which in turn is naturally 
isomorphic to $\cal H$$(H^*)$ (see \cite{lu1}), so that $D^{\omega }(H)^*_R$  
is indeed a generalization of the usual Heisenberg double of $H^*$.  

\end{document}